\documentclass{book}
\usepackage{amsmath,amsthm,amssymb}
\newtheorem{theorem}{Theorem}
\newtheorem{lemma}{Lemma}
\newtheorem{proposition}{Proposition}

\begin{document}

\chapter{On Landau's Function $g(n)$}
\label{ch18:chap18}

{Jean-Louis Nicolas}

\medskip

\noindent
Universit\'e de Lyon, CNRS, Universit\'e  Lyon 1,  
Institut Camille Jordan, Math\'ematiques, 
43 Bd. du 11 Novembre 1918,  
F-69622 Villeurbanne Cedex, France  

\noindent
\tt{jlnicola@in2p3.fr}   

\noindent
\tt{http://math.univ-lyon1.fr/$\sim$nicolas/}
\rm

\section{Introduction}
\label{ch18:sec18.1}

Let $S_n$ be the symmetric group of $n$ letters. Landau considered
the function $g(n)$ defined as the maximal order of an element of
$S_n$; Landau observed that (cf. \cite{ch18:bib5})
\begin{equation}
g(n)=\max\ {\rm lcm}(m_1,\ldots,m_k) \label{ch18:eqn1.1}
\end{equation}
where the maximum is taken on all the partitions $n = m_1 + m_2+
\ldots+m_k$ of $n$ and proved that, when $n$ tends to infinity
\begin{equation}
\log\ g(n)\sim\sqrt{n \log n}. \label{ch18:eqn1.2}
\end{equation} 
More precise asymptotic estimates have been given in
\cite{ch18:bib16,ch18:bib25,ch18:bib7}. In \cite{ch18:bib25} 
and  \cite{ch18:bib7} one also can find asymptotic estimates 
for the number of prime factors of $g(n)$. In \cite{ch18:bib26} 
and  \cite{ch18:bib29}, the largest prime factor $P^+(g(n))$ 
of $g(n)$ is investigated. In \cite{ch18:bib6} 
and  \cite{ch18:bib8}, effective upper and lower bounds 
of $g(n)$ are given. In \cite{ch18:bib28}, it is proved that
$\lim_{n\to\infty}g(n+1)/g(n)=1$. An algorithm able to calculate 
$g(n)$ up to $10^{15}$ is given in \cite{ch18:bib27} (see also
\cite{ch18:bib24}). The sequence of distinct values of $g(n)$ is 
entry A002809 of \cite{ch18:bib23}.  A nice survey paper was written by
W. Miller in 1987 (cf. \cite{ch18:bib9}).

My very first mathematical paper \cite{ch18:bib11}  was about
Landau's function, and the main result was that $g(n)$, which is
obviously non decreasing, is constant on arbitrarily long intervals
(cf. also~\cite{ch18:bib12}).
First time I met A.~Schinzel in Paris in May 1967. He told me that
he was interested in my results, but that P.~Erd\H{o}s would be more
interested than himself. Then I wrote my first letter to Paul with a
copy of my work. I received an answer dated of June 12 1967 saying
'' I sometimes thought about $g(n)$ but my results were very much less
complete than yours''. Afterwards, I met my advisor, the late
Professor Pisot, who, in view of this letter, told me that my work
was good  for a thesis.

The main idea of my work about $g(n)$ was to use the tools
introduced by S.~Ramanujan to study highly composite numbers (cf.
\cite{ch18:bib14,ch18:bib22}). P.~Erd\H{o}s was very well aware of this paper
of Ramanujan (cf. \cite{ch18:bib1,ch18:bib2,ch18:bib3,ch18:bib21}) 
as well as of the symmetric group and the order of
its elements, (cf. \cite{ch18:bib4}) and I think that he enjoyed the
connection between these two areas of mathematics. Anyway, since
these first letters, we had many occasions to discuss Landau's
function.

Let us define $n_1 = 1,n_2 = 2,n_3 = 3,n_4 = 4,n_5 = 5,n_6 = 7$, etc
$\ldots, n_k$ (see a table of $g(n)$ in \cite[p. 187]{ch18:bib12}), such that
\begin{equation}
g(n_k) > g(n_k -1). \label{ch18:eqn1.3}
\end{equation}
The above mentioned result can be read:
\begin{equation}
\overline{\lim}\,(n_{k+1} - n_k )= +\infty. \label{ch18:eqn1.4}
\end{equation}
Here, I shall prove the following result:
\begin{theorem}\label{ch18:them1.1}
\begin{equation}
\underline{\lim}\,( n_{k+1} - n_k) < +\infty. \label{ch18:eqn1.5}
\end{equation}
\end{theorem}
Let us set $p_1 = 2,p_2 = 3,p_3 = 5, \ldots$, $p_k=$ the $k$-th prime. It is
easy to deduce Theorem\,\ref{ch18:them1.1} from the twin prime
conjecture (i.e.\ $\underline{\lim}\,( p_{k+1}- p_k) = 2$) or even from
the weaker conjecture $\underline{\lim}\,( p_{k+1}- p_k) <+\infty $.
(cf.\ $\S$\ref{ch18:sec18.4} below). But I shall prove Theorem\,\ref{ch18:them1.1}
independently of these deep conjectures. Moreover I shall explain
below why it is reasonable to conjecture that the mean value of
$n_{k+1}- n_k$ is 2; in other terms one may conjecture that
\begin{equation}
n_k\sim 2k \label{ch18:eqn1.6}
\end{equation}
and that $n_{k+1}- n_k = 2$ has infinitely many solutions. Due to a
parity phenomenon, $n_{k+1}-n_k$ seems to be much more often even
than odd; nevertheless, I conjecture that:
\begin{equation}
\underline{\lim}\,( n_{k+1}-n_k) = 1. \label{ch18:eqn1.7}
\end{equation}

The steps of the proof of Theorem\,\ref{ch18:them1.1} are first to
construct the set $G$ of values of $g(n)$ corresponding to the so
called superior highly composite numbers introduced by S. Ramanujan,
and then, when $g(n) \in G$, to build the table of $g(n +d)$ when
$d$ is small. This will be done in \S \ref{ch18:sec18.4} 
and \S \ref{ch18:sec18.5}. Such values of
$g(n +d)$ will be linked with the number of distinct differences of
the form $P-Q$ where $P$ and $Q$ are primes satisfying
$x-x^{\alpha}\leq Q \leq x < P \leq x+x^{\alpha}$, where $x$ goes to
infinity and $0 < \alpha < 1$. Our guess is that these differences
$P - Q$ represent almost all even numbers between 0 and
$2x^{\alpha}$, but we shall only prove in \S \ref{ch18:sec18.3} 
that the number of
these differences is of the order of magnitude of $x^{\alpha}$,
under certain strong hypothesis on $x$ and $\alpha$, and for that a
result due to Selberg about the primes between $x$ and $x
+x^{\alpha}$ will be needed (cf. \S \ref{ch18:sec18.2}).

To support conjecture (\ref{ch18:eqn1.6}), I think that what has
been done here with $g(n) \in G$ can also be done for many more
values of $g(n)$, but, unfortunately, even assuming strong
hypotheses, I do not see for the moment how to manage it.

I thank very much E. Fouvry who gave me the proof of
Proposition\,\ref{ch18:prop3.1}.

\subsection{\bf{Notation}}
\label{ch18:subsec18.1.1}

$p$ will denote a generic prime, $p_k$ the $k$-th prime; $P,
Q,P_i,Q_j$ will also denote primes. As usual $\pi(x) =
\sum\nolimits_{p\leq x}1$ is the number of primes up to $x$.

$|S|$ will denote the number of elements of the set S. The sequence
$n_k$ is defined by (\ref{ch18:eqn1.3}).

\section{About the distribution of primes}
\label{ch18:sec18.2}

\begin{proposition}\label{ch18:prop2.1}
{\it Let us define} $\pi(x) = \sum\nolimits_{p\leq x} 1$, {\it and
let} $\alpha$ {\it be such that} $\frac{1}{6}<\alpha< 1$, {\it and}
$\varepsilon  > 0$. {\it When} $\xi$ \textit{goes to infinity, and}
$\xi^{\prime}= \xi+\xi/\log \xi$, \textit{then for all x in the
interval} $[\xi,\xi^{\prime}]$ \textit{but a subset of measure}
$O((\xi^{\prime}- \xi)/ \log^{3} \xi)$ {\it we have}:

\begin{equation}
\left|\pi(x+x^{\alpha})-\pi(x)-\frac{x^{\alpha}}{\log x}\right|\leq \varepsilon
\frac{x^{\alpha}}{\log x} \label{ch18:eqn2.1}
\end{equation}

\begin{equation}
\left|\pi(x)-\pi(x-x^{\alpha})-\frac{x^{\alpha}}{\log x}\right|\leq\varepsilon
\frac{x^{\alpha}}{\log x} \label{ch18:eqn2.2}
\end{equation}

\begin{equation}
\left|\frac{x}{\log x}-\frac{Q^{k}-Q^{k-1}}{\log Q}\right|\geq
\frac{\sqrt{x}}{\log^{4}x}\ for\ all\ primes\ Q,\ and\ k\geq 2.
\label{ch18:eqn2.3}
\end{equation}
\end{proposition}

{\it Proof}. This proposition is an easy extension of a result of Selberg
(cf. \cite{ch18:bib15}) who proved that (\ref{ch18:eqn2.1}) holds
for most $x$ in $(\xi,\xi^{\prime})$. In \cite{ch18:bib13}, I gave a
first extension of Selberg's result by proving that
(\ref{ch18:eqn2.1}) and (\ref{ch18:eqn2.2}) hold simultaneously for
all $x$ in $(\xi,\xi^{\prime})$ but for a subset of measure
$O((\xi^{\prime} - \xi)/ \log^{3}\xi)$. So, it suffices to prove
that the measure of the set of values of $x$ in $(\xi,\xi^{\prime})$
for which (\ref{ch18:eqn2.3}) does not hold is $O((\xi^{\prime} -
\xi)/ \log^{3}\xi)$.

We first count the number of primes $Q$ such that for one $k$ we
have:
\begin{equation}
\frac{\xi}{\log \xi}\leq\frac{Q^{k}-Q^{k-1}}{\log Q}\leq
\frac{\xi^{\prime}}{\log \xi^{\prime}}. \label{ch18:eqn2.4}
\end{equation}
If $Q$ satisfies (\ref{ch18:eqn2.4}), then $k \leq\frac{\log
\xi^{\prime}}{\log 2}$ for $\xi^{\prime}$ large enough. Further, for
$k$ fixed, (\ref{ch18:eqn2.4}) implies that $Q\leq
(\xi^{\prime})^{1/k}$, and the total number of solutions of
(\ref{ch18:eqn2.4}) is
\[
\leq\sum\limits_{k=2}^{\log \xi^{\prime}/\log
2}(\xi^{\prime})^{1/k}=O(\sqrt{\xi^{\prime}})=O(\sqrt{\xi}).
\]
With a more careful estimation, this upper bound could be improved, 
but this crude result is enough for our purpose. Now,
for all values of $y =\frac{Q^{k}-Q^{k-1}}{\log Q} $ satisfying
(\ref{ch18:eqn2.4}), we cross out the interval
$\left(y-\frac{\sqrt{\xi^{\prime}}}{\log^{4}\xi^{\prime}},y+\frac{\sqrt{\xi^{\prime}}}{\log^{4}\xi^{\prime}}\right)$.
We also cross out this interval whenever $y = \frac{\xi}{\log \xi}$
and $y = \frac{\xi^{\prime}}{\log \xi^{\prime}}$. The total sum of
the lengths of the crossed out intervals is $O
\left(\frac{\xi}{\log^{4} \xi}\right)$, which is smaller than the
length of the interval $\left(\frac{\xi}{\log
\xi},\frac{\xi^{\prime}}{\log \xi^{\prime}}\right)$ and if
$\frac{x}{\log x}$ does not fall into one of these forbidden
intervals, (\ref{ch18:eqn2.3}) will certainly hold. Since the
derivative of the function $\varphi(x) = x/ \log x$ is
$\varphi^{\prime}(x) = \frac{1}{\log x}-\frac{1}{\log^{2} x}$ and
satisfies $\varphi^{\prime}(x)\sim\frac{1}{\log \xi}$ for all $x \in
(\xi,\xi^{\prime})$, the measure of the set of values of $x \in
(\xi,\xi^{\prime})$ such that $\varphi(x)$ falls into one of the
above forbidden intervals is, by the mean value theorem
$O\left(\frac{\xi}{\log^{3}\xi}\right)$, and the proof of
Proposition\,\ref{ch18:prop2.1} is completed.

\section{About the differences between primes}
\label{ch18:sec18.3}

\begin{proposition}\label{ch18:prop3.1}
\textit{Suppose that there exists} $\alpha, 0 < \alpha< 1$, {\it
and} $x$ \textit{large enough such that the inequalities}

\begin{equation}
\pi(x+x^{\alpha})-\pi(x)\geq(1-\varepsilon)x^{\alpha}/\log x
\label{ch18:eqn3.1}
\end{equation}

\begin{equation}
\pi(x)-\pi(x-x^{\alpha})\geq(1-\varepsilon)x^{\alpha}/\log x
\label{ch18:eqn3.2}
\end{equation}

\textit{hold. Then the set}
\[
E = E(x, \alpha) = \{P- Q;P,Q\ primes,\ x -x^{\alpha}< Q\leq x <
P\leq x+ x^{\alpha}\}
\]
{\it satisfies:}
\[
|E|\geq C_2 x^{\alpha}
\]
{\it where} $C_2 = C_1 \alpha^{4}(1- \varepsilon)^{4}$ {\it and}
$C_1$ \textit{is an absolute constant} $(C_1 = 0.00164\ works)$.
\end{proposition}

{\it Proof}. The proof is a classical application of the sieve method that
Paul Erd\H{o}s enjoys very much. Let us set, for $d \leq
2x^{\alpha}$,
\[
r(d) =|\{(P,Q);x - x^{\alpha}< Q\leq x < P\leq x + x^{\alpha},P - Q
= d\}|.
\]
Clearly we have
\begin{equation}
|E|=\sum\limits_{\substack{0<d\leq2x^{\alpha}\\r(d)\neq 0}} 1
\label{ch18:eqn3.3}
\end{equation}
and
\begin{equation}
\sum\limits_{0<d\leq2x^{\alpha}}r(d)=(\pi(x+x^{\alpha})-\pi(x))(\pi(x)-\pi(x-x^{\alpha}))\geq(1-\varepsilon)^{2}x^{2\alpha}/\log^{2}
x. \label{ch18:eqn3.4}
\end{equation}
Now to get an upper bound for $r(d)$, we sift the set
\[
A = \{n; x - x^{\alpha}< n \leq x\}
\]
with the primes $p \leq z$. If $p$ divides $d$, we cross out the
$n^{\prime}s$ satisfying $n \equiv\ 0 \pmod{p}$, and if $p$ does not
divide $d$, the $n^{\prime}s$ satisfying
\[
n \equiv 0 \pmod{p}\quad  {\rm or}\quad n \equiv -d \pmod{p}
\]
so that we set for $p \leq z$:
$$
w(p)=
\begin{cases}
1& {\rm if}\ p\ {\rm divides}\ d\\
2& {\rm if}\ p\ {\rm does\ not\ divide}\ d.
\end{cases}
$$
By applying the large sieve (cf. \cite[Corollary 1]{ch18:bib10}), we have
\[
r(d)\leq\frac{|A|}{L(z)}
\]
with
\[
L(z)=\sum\limits_{n\leq
z}\left(1+\frac{3}{2}n|A|^{-1}z\right)^{-1}\mu(n)^{2}\left(\mathop
{\prod}\limits_{p|n}\frac{w(p)}{p-w(p)}\right)
\]
($\mu$ is the M\"{o}bius function), and with the choice $z =
(\frac{2}{3}|A|)^{1/2}$, it is proved in \cite{ch18:bib17} that
\[
\frac{|A|}{L(z)}\leq 16\mathop{\prod}\limits_{p\geq
3}\left(1-\frac{1}{(p-1)^{2}}\right)\frac{|A|}{\log^{2}(|A|)}\mathop{\prod}\limits_{\substack{p|d\\p>2}}\frac{p-1}{p-2}.
\]
The value of the above infinite product is $0.6602 \ldots < 2/3$. We
set $f(d) = \prod_{\substack{p|d\\p>2}}\frac{p-1}{p-2}$, and we
observe that $|A|\geq x^{\alpha}- 1$, so that for $x$ large enough
\begin{equation}
r(d)\leq\frac{32}{3\alpha^{2}}\frac{|A|}{\log^{2}x}f(d).
\label{ch18:eqn3.5}
\end{equation}
Now, for the next step, we shall need an upper bound for $\sum_{n\leq
x}f^{2}(n)$. By using the convolution method and defining
\[
h(n) = \sum\limits_{a|n}\mu(a)f^{2}(n/a)
\]
one gets $h(2) = h(2^2) = h(2^3) =\ldots = 0$ and, for $p \geq 3$,
$h(p) =\frac{2p-3}{(p-2)^{2}}$, $h(p^{2}) = h(p^{3}) = \ldots= 0$, so that
\begin{eqnarray}
\sum\nolimits_{n\leq x}f^{2}(n)&=&\sum\nolimits_{n\leq x}\sum\nolimits_{a|n}h(a)=\sum\nolimits_{a\leq x}h(a)\left\lfloor{\frac{x}{a}}\right\rfloor\nonumber\\
&\leq&x\sum\nolimits_{a=1}^{\infty}\frac{h(a)}{a}=x\prod\nolimits_{p\geq3}\left(1+\frac{2p-3}{p(p-2)^{2}}\right)\label{ch18:eqn3.6}\\
&=&2.63985 \ldots x\leq\frac{8}{3}x.\nonumber
\end{eqnarray}
From (\ref{ch18:eqn3.4}) and (\ref{ch18:eqn3.5}), one can deduce
\[
\frac{(1-\varepsilon)^{2}x^{2\alpha}}{\log^{2}x}\leq\sum\limits_{\substack{0<d\leq2x^{\alpha}\\r(d)\neq
0}}r(d)\leq\frac{32}{3\alpha^{2}}\frac{|A|}{\log^{2}x}\sum\limits_{\substack{0<d\leq2x^{\alpha}\\r(d)\neq
0}}f(d).
\]
which implies
\[
\sum\limits_{\substack{0<d\leq2x^{\alpha}\\r(d)\neq
0}}f(d)\geq\frac{3\alpha^{2}x^{2\alpha}(1-\varepsilon)^{2}}{32|A|}\,\cdot
\]
By Cauchy-Schwarz's inequality, one has
\[
\left(\sum\limits_{\substack{0<d\leq2x^{\alpha}\\r(d)\neq
0}}1\right)\left(\sum\limits_{\substack{0<d\leq2x^{\alpha}\\r(d)\neq
0}}f^{2}(d)\right)\geq\frac{9\alpha^{4}x^{4\alpha}(1-\varepsilon)^{4}}{1024|A|^{2}}
\]
and, by (\ref{ch18:eqn3.3}) and (\ref{ch18:eqn3.6})
\[
|E|\geq\frac{9\alpha^{4}x^{4\alpha}(1-\varepsilon)^{4}}
{1024|A|^{2}}\left/\frac 83 (2x^\alpha)
\right. =\frac{27}{16384} \frac{x^{3\alpha}(1-\varepsilon)^{4}}{|A|^2}\cdot
\]
Since $|A|\leq x^{\alpha}+ 1$, and $x$ has been supposed large
enough, proposition \ref{ch18:prop3.1} is proved.

\section{Some properties of $g(n)$}
\label{ch18:sec18.4}

Here, we recall some known properties of $g(n)$ which can be found
for instance in \cite{ch18:bib12}. Let us define the arithmetic
function $\ell$ in the following way: $\ell$ is additive, and, if $p$ is a
prime and $k \geq 1$, then $\ell(p^{k}) = p^{k}$. It is not difficult
to deduce from (1.1) (cf. \cite{ch18:bib9} or \cite{ch18:bib12})
that
\begin{equation}
g(n) = \max_{\ell(M)\leq n}M. \label{ch18:eqn4.1}
\end{equation}
Now the relation (cf. \cite{ch18:bib12}, p. 139)
\begin{equation}
M\in g(\mathbb{N}) \quad\Longleftrightarrow \quad (M^{\prime} > 
M \;\;\Longrightarrow\;\;
\ell(M^{\prime}) > \ell(M)) \label{ch18:eqn4.2}
\end{equation}
easily follows from (\ref{ch18:eqn4.1}), and shows that the values
of the Landau function $g$ are the ''champions'' for the small
values of $\ell$. So the methods introduced by Ramanujan (cf.
\cite{ch18:bib14}) to study highly composite numbers can also be
used for $g(n)$. Indeed $M$ is highly composite, if it is a
''champion'' for the divisor function $d$, that is to say if
\[
M^{\prime} < M \quad\Longrightarrow\quad d(M^{\prime}) < d(M).
\]
Corresponding to the so-called superior highly composite numbers,
one introduces the set $G : N\in G$ if there exists $\rho > 0$ such
that
\begin{equation}
\forall M \geq 1,\quad \ell(M) - \rho \log M \geq \ell(N) - \rho \log N.
\label{ch18:eqn4.3}
\end{equation}
(\ref{ch18:eqn4.2}) and (\ref{ch18:eqn4.3}) easily imply that $G
\subset g(\mathbb{N})$. Moreover, if $\rho > 2/ \log 2$, let us
define $x > 4$ such that $\rho = x/\log x$ and
\begin{equation}
N_{\rho} = \mathop{\prod}\limits_{p\leq
x}p^{\alpha_{p}} = \mathop{\prod}\limits_{p}p^{\alpha_{p}}
\label{ch18:eqn4.4}
\end{equation}
with
\[
\alpha_{p}=
\begin{cases}
0& {\rm if} \quad p> x\\
1&{\rm if} \quad \frac{p}{\log p}\leq\rho<\frac{p^{2}-p}{\log p}\\
k\geq2&{\rm if}\quad\frac{p^{k}-p^{k-1}}{\log
p}\leq\rho<\frac{p^{k-1}-p^{k}}{\log p}
\end{cases}
\]
then $N_{\rho} \in G$. With the above definition, since $x\geq 4$,
it is not difficult to show that (cf. \cite[(5)]{ch18:bib7})
\begin{equation}
p^{\alpha_{p}}\leq x \label{ch18:eqn4.5}
\end{equation}
holds for $p \leq x$, whence $N_{\rho}$ is a divisor of the l.c.m.
of the integers $\leq x$. Here we can prove

\begin{proposition}\label{ch18:prop4.1}
\textit{For every prime p, there exists n such that the largest
prime factor of g(n) is equal to p}.
\end{proposition}

{\it Proof}. We have $g(2) = 2, g(3) = 3$. If $p \geq 5$, let us choose
$\rho = p/ \log p > 2/\log 2$. $N_{\rho}$ defined by
(\ref{ch18:eqn4.4}) belongs to $G \subset g(\mathbb{N})$, and its
largest prime factor is $p$, which proves Proposition
\eqref{ch18:prop4.1}.

\medskip

From Proposition\,\ref{ch18:prop4.1}, it is easy to deduce a proof
of Theorem\,\ref{ch18:them1.1}, under the twin prime conjecture. Let
$P = p + 2$ be twin primes, and $n$ such that the largest prime
factor of $g(n)$ is $p$. The sequence $n_k$ being defined by
(\ref{ch18:eqn1.3}), we define $k$ in terms of $n$ by $n_k \leq n <
n_{k+1}$, so that $g(n_k) = g(n)$ has its largest prime factor equal
to $p$. Now, from (\ref{ch18:eqn4.1}) and (\ref{ch18:eqn4.2}),
\[
\ell(g(n_k)) = n_k
\]
and $g(n_k + 2) > g(n_k)$ since $M = \frac{P}{p} g(n_k)$ satisfies $M 
> g(n_k)$ and $\ell(M) =n_k + 2$. So $n_{k+1}\leq n_k + 2$, and
Theorem\,\ref{ch18:them1.1} is proved under this strong hypothesis.

\medskip

Let us introduce now the so-called benefit method. For a fixed $\rho
> 2/\log 2$, $N = N_{\rho}$ is defined by (\ref{ch18:eqn4.4}), and for any
integer $M$,
\[
M=\mathop{\prod}\limits_{p}p^{\beta_{p}},
\]
one defines the benefit of $M$:
\begin{equation}
{\rm ben}(M) = \ell(M) -\ell(N) - \rho \log M/N. \label{ch18:eqn4.6}
\end{equation}
Clearly, from (\ref{ch18:eqn4.3}), ${\rm ben}( M) \geq 0$ holds, and 
from the additivity of $\ell$ one has 
\begin{equation}
{\rm ben}(M) = \sum\limits_{p}\left(\ell(p^{\beta_{p}})-\ell(p^{\alpha_{p}}) -
\rho(\beta_{p} -\alpha_{p})\log p\right). \label{ch18:eqn4.7}
\end{equation}
In the above formula, let us observe that $\ell(p^{\beta}) = p^{\beta}$
if $\beta\geq 1$, but that $\ell(p^{\beta}) = 0 \neq p^{\beta}= 1$ if
$\beta=0$, and, due to the choice of $\alpha_{p}$ in
(\ref{ch18:eqn4.4}), that, in the sum \eqref{ch18:eqn4.7},
all the terms are non negative:
for all $p$ and for $\beta \geq 0$, we have
\begin{equation}
\ell(p^{\beta})-\ell(p^{\alpha_p})-\rho(\beta-\alpha_p)\log p\geq 0
\label{ch18:eqn4.8}
\end{equation}

Indeed, let us consider the set of points (0,0) and 
$(\beta,p^{\beta}/ \log p)$ for $\beta$
integer $\geq 1$. For all $p$, the
piecewise linear curve going through these points is convex, and for
a given $\rho,\alpha_{p}$ is chosen so that the straight line $L$ of
slope $\rho$ going through
$\left(\alpha_{p},\frac{p^{\alpha_{p}}}{\log p}\right)$ does not cut
that curve. The left-hand side of (\ref{ch18:eqn4.8}), (which is
${\rm ben}(Np^{\beta-\alpha_{p}})$) can be seen as the product of log
$p$ by the vertical distance of the point
$\left(\beta,\frac{p^{\beta}}{\log p}\right)$ to the straight
line $L$, and because of convexity, we shall have for all $p$,
\begin{equation}
{\rm ben}(Np^{t}) \geq t\ {\rm ben}(Np),\quad  t \geq 1 \label{ch18:eqn4.9}
\end{equation}
and for $p \leq x$,
\begin{equation}
{\rm ben}(N{p}{^{-t}})\geq t\ {\rm ben}(N{p}{^{-1}}),\quad 1\leq t\leq\alpha_p.
\label{ch18:eqn4.10}
\end{equation}

\section{Proof of Theorem\,\ref{ch18:them1.1}}
\label{ch18:sec18.5}

First the following proposition will be proved:

\begin{proposition}\label{ch18:prop5.1}
{\it Let} $\alpha< 1/2$, {\it and} $x$ \textit{large enough such
that (\ref{ch18:eqn2.3}) holds. Let us denote the primes surrounding
x by}:
\[
\ldots < Q_{j} < \ldots < Q_2 < Q_1 \leq x < P_1 < P_2 < \ldots <
P_i< \ldots
\]
\textit{Let us define} $\rho = x/ \log x, N = N_{\rho}$ {\it by}
(\ref{ch18:eqn4.4}), $n = \ell(N)$. {\it Then for} $n \leq m \leq n +
2x^{\alpha}, g(m)$ \textit{can be written}
\begin{equation}
g(m)=N\frac{P_{i_{1}}P_{i_{2}}\ldots
P_{i_{r}}}{Q_{j_{1}}Q_{j_{2}}\ldots Q_{j_{r}}} \label{ch18:eqn5.1}
\end{equation}
{\it with} $r\geq 0$ {\it and} $i_1<\ldots<i_r,j_1<\ldots<j_r,P_{i_{r}}\leq
x+4x^{\alpha},Q_{j_{r}}\geq x- 4x^{\alpha}.$
\end{proposition}

{\it Proof}. First, from (\ref{ch18:eqn4.1}), one has $\ell(g(m)) \leq m$, and
from (\ref{ch18:eqn4.6}) and (\ref{ch18:eqn4.1})
\begin{equation}
{\rm ben}(g(m)) = \ell(g(m)) -\ell(N) - \rho \log \frac{g(m)}{N} \leq m - n \leq
2x^{\alpha} \label{ch18:eqn5.2}
\end{equation}
for $n \leq m \leq 2x^{\alpha}$.
\medskip

Further, let $Q \leq x$ be a prime, and $k = \alpha_{Q} \geq 1$ the
exponent of $Q$ in the standard factorization of $N$. Let us suppose
that for a fixed $m, Q$ divides $g(m)$ with the exponent $\beta_{Q}
= k +t, t > 0$. Then, from (\ref{ch18:eqn4.7}), (\ref{ch18:eqn4.8}),
and (\ref{ch18:eqn4.9}), one gets
\begin{equation}
{\rm ben}(g(m)) \geq {\rm ben}(NQ^{t})\geq {\rm ben}(NQ)\label{ch18:eqn5.3}
\end{equation}
and
\begin{eqnarray*}
{\rm ben}(NQ)&=&Q^{k+1}-Q^{k}-\rho \log Q\\
&=&\log Q\left(\frac{Q^{k+1}-Q^{k}}{\log Q}-\rho\right).
\end{eqnarray*}
From (\ref{ch18:eqn4.4}), the above parenthesis is non negative, and
from (\ref{ch18:eqn2.3}), one gets:
\begin{equation}
{\rm ben}(NQ)\geq \log 2\frac{\sqrt{x}}{\log^{4}x}\cdot \label{ch18:eqn5.4}
\end{equation}
For $x$ large enough, there is a contradiction between
(\ref{ch18:eqn5.2}), (\ref{ch18:eqn5.3}) and (\ref{ch18:eqn5.4}),
and so, $\beta_{Q} \leq \alpha_{Q}$.
\medskip

Similarly, let us suppose $Q \leq x, k = \alpha_{Q} \geq 2$ and 
$\beta_{Q} = k - t, 1 \leq t \leq k$. One has, from
(\ref{ch18:eqn4.7}), (\ref{ch18:eqn4.8}) and (\ref{ch18:eqn4.10}),
\[
{\rm ben}(g(m)) \geq {\rm ben} (NQ^{-t}) \geq {\rm ben}(NQ^{-1})
\]
and
\begin{eqnarray*}
{\rm ben}(NQ^{-1})&=&Q^{k-1}-Q^{k}+\rho \log Q\\
&=&\log Q\left(\rho-\frac{Q^{k}-Q^{k-1}}{\log Q}\right)
\geq \log 2\frac{\sqrt{x}}{\log^{4}x}
\end{eqnarray*}
which contradicts (\ref{ch18:eqn5.2}), and so, for such a
$Q$, $\beta_{Q} = \alpha_{Q}$.
\medskip

Now, let us suppose $Q \leq x, \alpha_{Q} = 1$, and $\beta_{Q} = 0$ for
some $m, n \leq m \leq n + 2x^{\alpha}$. Then
\begin{equation*}
{\rm ben}(g(m)) \geq {\rm ben}(NQ^{-1}) = -Q +\rho \log Q = y(Q)
\end{equation*}
by setting $y(t)=\rho \log t-t$. From the concavity of $y(t)$ for $t >
0$, for $x\geq e^2$, we get
\begin{eqnarray*}
y(Q)\geq y(x) + (Q-x) y'(x) &=& (Q-x) \left( \frac \rho x-1\right)\\
&=& (x-Q)\left(1-\frac{1}{\log x}\right)\geq \frac 12 (x-Q)
\end{eqnarray*}
and  so,
\[
{\rm ben}(g(m))\geq\frac{1}{2}(x-Q)
\]
which, from (\ref{ch18:eqn5.2}) yields
\[
x-Q\leq 4x^{\alpha}.
\]
In conclusion, the only prime factors allowed
in the denominator of $\frac{g(m)}{N}$ are the $Q^{\prime}s$, with
$x - 4x^{\alpha}\leq Q \leq x$ , and $\alpha_{Q} = 1$.
\medskip

What about the numerator? Let $P > x$ be a prime number and suppose
that $P^{t}$ divides $g(m)$ with $t \geq 2$. Then, from
(\ref{ch18:eqn4.9}) and (\ref{ch18:eqn4.6}),
\[
{\rm ben}(Np^{t}) \geq {\rm ben}(Np^{2}) = P^{2} - 2\rho \log P.
\]
But the function $t \mapsto t^{2} - 2 \rho \log t$ is increasing for
$t \geq \sqrt{\rho}$, so that,
\[
{\rm ben}(NP^{t})\geq x^{2}-2x>2x^{\alpha}
\]
for $x$ large enough, which contradicts (\ref{ch18:eqn5.2}). The
only possibility is that $P$ divides $g(m)$ with exponent 1. In that
case, from the convexity of the function $z(t)=t-\rho \log t$,
inequality (\ref{ch18:eqn4.9}) yields
\begin{eqnarray*}
{\rm ben}(g(m)) \;\geq \;
{\rm ben}(NP) \hspace{-2mm} &=& \hspace{-2mm} 
 z(P) \; \geq \; z(x) + (P-x)z'(x)\\
&=& \hspace{-2mm}
(P - x) \left(1-\frac{1}{\log x}\right) \; \geq \;\frac{1}{2}(P - x)
\end{eqnarray*}
for $x\geq e^2$, which, with (\ref{ch18:eqn5.2}),  implies
\[
P-x\leq 4x^{\alpha}.
\]

Up to now, we have shown that
\[
g(m)=N\frac{P_{i_{1}}\ldots P_{i_{r}}}{Q_{j_{1}}\ldots Q_{j_{s}}}
\]
with $P_{i_{r}}\leq x + 4x^{\alpha},Q_{j_{s}}\geq x - 4x^{\alpha}$.
It remains to show that $r = s$. First, since $n \leq m \leq n +
2x^{\alpha}$, and $N$ belongs to $G$, we have from
(\ref{ch18:eqn4.1}) and (\ref{ch18:eqn4.2})
\begin{equation}
n \leq \ell(g(m)) \leq n + 2x^{\alpha}. \label{ch18:eqn5.5}
\end{equation}
Further,
\[
\ell(g(m))-n=\sum\limits_{t=1}^{r}P_{i_{t}}-\sum\limits_{t=1}^{s}Q_{j_{t}}
\]
and since $r \leq 4x^{\alpha}$, and $s \leq4x^{\alpha}$,
\begin{eqnarray*}
\ell(g(m))-n&\leq&r(x+4x^{\alpha})-s(x-4x^{\alpha})\\
&\leq&(r-s)x+32x^{2\alpha}.
\end{eqnarray*}
From \eqref{ch18:eqn5.5}, $\ell(g(m))-n \geq 0$ holds and
as $\alpha< 1/2$, this implies that $r \geq s$
for $x$ large enough. Similarly,
\[
\ell(g(m))- n \geq (r - s)x,
\]
so, from  \eqref{ch18:eqn5.5}, $(r - s)x$ must be 
$\leq 2x^{\alpha}$,  which, for $x$ large
enough, implies $r \leq s$; finally $r = s$, and the proof
of Proposition\,\ref{ch18:prop5.1} is completed.

\begin{lemma}\label{ch18:lem5.1}
\textit{Let x be a positive real number}, $a_1, a_2, \ldots, a_k,
b_1, b_2, \ldots,b_k$ \textit{be real number such that}
\[
b_k \leq b_{k-1} \leq \ldots\leq b_1 \leq x < a_1 \leq a_2 \leq
\ldots \leq a_k
\]
{\it and} $\Delta$ {\it be defined by} $\Delta =\sum_{i=1}^{k}(a_i-
b_i)$. \textit{Then the following inequalities}
\[
\frac{x+\Delta}{x}\leq\mathop{\prod}\limits^{k}_{i=1}\frac{a_i}{b_i}\leq
\exp\left(\frac{\Delta}{x}\right)
\]
hold.
\end{lemma}

{\it Proof}. It is easy, and can be found in \cite{ch18:bib12},
p. 159.

\medskip

Now it is time to prove Theorem\,\ref{ch18:them1.1}. With the
notation and hypothesis of Proposition\,\ref{ch18:prop5.1}, let us
denote by $B$ the set of integers $M$ of the form
\[
M=N\frac{P_{i_{1}}P_{i_{2}}\ldots
P_{i_{r}}}{Q_{j_{1}}Q_{j_{2}}\ldots Q_{j_{r}}}
\]
satisfying
\[
\ell(M) -\ell(N) = \sum\limits_{t=1}^{r}(P_{i_{t}} - {Q_{j_{t}}})\leq
2x^{\alpha}.
\]
From Proposition\,\ref{ch18:prop5.1}, for $n \leq m \leq
2x^{\alpha}, g(m) \in B$, and thus, from (\ref{ch18:eqn4.1}),
\begin{equation}
g(m) = \mathop {\max\nolimits _{\ell(M) \le m}M.}\limits_{M \in B}
\label{ch18:eqn5.6}
\end{equation}
Further, for $0\leq d \leq 2x^{\alpha}$, define
\[
B_{d} = \{M \in B ;\ell(M) -\ell(N) = d\}.
\]
I claim that, if $d < d^{\prime}$ (which implies $d \leq d^{\prime}-
2)$, any element of $B_d$ is smaller than any element of
$B_{d^{\prime}}$. Indeed, let $M\in B_d$, and $M^{\prime}\in
B_{d^{\prime}}$. From  Lemma \ref{ch18:lem5.1}, one has
\[
\frac{M}{N}\leq \exp\left(\frac{d}{x}\right)\quad {\rm and}\quad
\frac{M^{\prime}}{N}\geq\frac{x+d^{\prime}}{x}\geq\frac{x+d+2}{x}\cdot
\]
Since $d < 2x^{\alpha}< x$, and $e^{t} \leq \frac{1}{1-t}$ for $0\leq t <
1$, one gets
\[
\frac{M}{N}\leq\frac{1}{1-d/x}=\frac{x}{x-d}\cdot
\]
This last quantity is smaller than $\frac{x+d+2}{x}$ if $(d + 1)^{2}
< 2x + 1$, which is true for $x$ large enough, because 
$d \leq 2x^{\alpha}$ and $\alpha <  1/2$.

From the preceding claim, and from (\ref{ch18:eqn5.6}), it follows
that, if $B_d$ is non empty, then
\[
g(n + d) = \max B_{d}.
\]
Further, since $N\in G$, we know that $n = \ell(N)$ belongs to the
sequence $(n_k)$ where $g$ is increasing, and so, $n = n_{k_{0}}$.
If $0 < d_1 < d_2 < \ldots < d_s\leq 2x^{\alpha}$ denote the values
of $d$ for which $B_d$ is non empty, then one has
\begin{equation}
n_{k_{0+i}}=n+d_i, 1\leq i\leq s. \label{ch18:eqn5.7}
\end{equation}

Suppose now that $\alpha < 1/2$ and $x$ have been chosen in such a
way that (\ref{ch18:eqn3.1}) and (\ref{ch18:eqn3.2}) hold. With the
notation of Proposition\,\ref{ch18:prop3.1}, the set $E(x , \alpha)$
is certainly included in the set $\{d_1 ,d_2,\ldots, d_s\}$, and
from Proposition\,\ref{ch18:prop3.1},
\begin{equation}
s\geq C_2 x^{\alpha} \label{ch18:eqn5.8}
\end{equation}
which implies that for at least one $i, d_{i+1}- d_{i}\leq
\frac{2}{C_2}$, and thus
\[
n_{k_{0}+i+1}-n_{k_{0}+i}\leq\frac{2}{C_2}.
\]
Finally, for $\frac{1}{6}< \alpha< \frac{1}{2}$,
Proposition\,\ref{ch18:prop2.1} allows us to choose $x$ as wished,
and thus, the proof of Theorem\,\ref{ch18:them1.1} is completed.
With $\varepsilon$ very small, and $\alpha$ close to $1/2$, the
values of $C_1$ and $C_2$ given in Proposition\,\ref{ch18:prop3.1}
yield that for infinitely many $k^{\prime}s$,
\[
n_{k+1}-n_k\leq20000.
\]

To count how many such differences we get, we define
$$\gamma(n) = { \rm Card}\{m\leq n ; g(m) > g(m-1)\}.$$ 
Therefore, with the notation (\ref{ch18:eqn1.3}), we have
$n_{\gamma_{(n)}}= n$. 

In \cite[162--164]{ch18:bib12}, it is proved that
\[
n^{1-\tau/2}\ll \gamma(n)\leq n-c\frac{n^{3/4}}{\sqrt{\log n}}
\]
where $\tau$ is such that the sequence of consecutive primes satisfies
$p_{i+1}- p_i \ll p_{i}^{\tau}$.
Without any hypothesis, the best known $\tau$ is $> 1/2$.

\begin{proposition}\label{ch18:prop5.2}
{\it We have} $\gamma(n) \geq n^{3/4-\varepsilon}$ {\it for all} $
\varepsilon> 0$, \textit{and n large enough}.
\end{proposition}

{\it Proof}. With the definition of $\gamma(n)$, (\ref{ch18:eqn5.7}) and
(\ref{ch18:eqn5.8}) give
\begin{equation}
\gamma(n+2x^{\alpha})-\gamma(n)\geq s\gg x^{\alpha}
\label{ch18:eqn5.9}
\end{equation}
whenever $n = \ell(N), N = N_{\rho}, \rho = x/ \log x$, and $x$
satisfies Proposition\,\ref{ch18:prop2.1}. But, from
(\ref{ch18:eqn4.4}), two close enough distinct values of $x$ can
yield the same $N$.

I now claim that, with the notation of
Proposition\,\ref{ch18:prop2.1}, the number of primes $p_i$ between
$\xi$ and $\xi^{\prime}$ such that there is at least one $x \in
[p_i, p_{i+1})$ satisfying (\ref{ch18:eqn2.1}), (\ref{ch18:eqn2.2})
and (\ref{ch18:eqn2.3}) is bigger than
$\frac{1}{2}(\pi(\xi^{\prime})-\pi(\xi))$. Indeed, for each $i$ for
which $[p_i,p_{i+1})$ does not contain any such $x$, we get a
measure $p_{i+1}-p_{i}\geq 2$, and if there are more than
$\frac{1}{2}(\pi(\xi^{\prime})-\pi(\xi))$ such $i^{\prime}s$, the
total measure will be greater than
$\pi(\xi^{\prime})-\pi(\xi)\sim\xi/\log^{2}\xi$, which contradicts
Proposition\,\ref{ch18:prop2.1}.

From the above claim, there will be at least
$\frac{1}{2}(\pi(\xi^{\prime})-\pi(\xi))$ distinct $N^{\prime}s$,
with $N = N_{\rho},\rho = x/\log x$, and $\xi\leq x \leq
\xi^{\prime}$. Moreover, for two such distinct $N$, say $N^{\prime}<
N$'', we have from (\ref{ch18:eqn4.4}), 
$\ell(N'')-\ell(N^\prime) \geq \xi$.

Let $N^{(1)}$ and $N^{(0)}$ the biggest and the smallest of these
$N^{\prime}s$, and $n^{(1)} = \ell(N^{(1)}), n^{(0)} = \ell(N^{(0)})$,
then from (\ref{ch18:eqn5.9}),
\begin{equation}
\gamma(n^{(1)})\geq\gamma(n^{(1)})-\gamma(n^{(0)})
\geq \frac 12 \left(\pi(\xi')-\pi(\xi)\right)\xi^\alpha
\gg \frac{\xi^{1+\alpha}}{\log^{2}\xi}. \label{ch18:eqn5.10}
\end{equation}
But from (\ref{ch18:eqn4.4}) and (\ref{ch18:eqn4.5}), $x\sim \log
N_{\rho}$, and from (\ref{ch18:eqn1.2}),
\[
x\sim\log N_{\rho}\sim\sqrt{n \log n}\quad {\rm with}\quad n=\ell(N_{p})
\]
so
\[
\xi\sim\sqrt{n^{(1)}\log n^{(1)}}
\]
and since $\alpha$ can be choosen in (\ref{ch18:eqn5.10}) as close
as wished of $1/2$, this completes the proof of
Proposition\,\ref{ch18:prop5.2}.


\begin{thebibliography}{17}

\bibitem{ch18:bib1} L. Alaoglu, P. Erd\H{o}s, 
``On highly composite and similar numbers'', 
Trans. Amer. Math. Soc. 56, 1944, 448--469.

\bibitem{ch18:bib27} M. Del\'eglise,  J.-L. Nicolas, P. Zimmermann, 
``Landau's function for one million billions'', 
J. de Th\'eorie des Nombres de Bordeaux, 20, 2008, 625--671.

\bibitem{ch18:bib29} M. Del\'eglise, J.-L. Nicolas, 
``Le plus grand facteur premier de la fonction de Landau'', 
Ramanujan J., 27, 2012, 109--145.

\bibitem{ch18:bib2} P. Erd\H{o}s, 
``On highly composited numbers'', 
J. London Math. Soc., 19, 1944, 130--133.

\bibitem{ch18:bib21} P. Erd\H{o}s, ``Ramanujan and I'', Number Theory,
  Madras 1987, Editor : K. Alladi, Lecture Notes in Mathematics n$^o$
  1395, Springer-Verlag, 1989.

\bibitem{ch18:bib3} P. Erd\H{o}s, J.-L. Nicolas, 
``R\'{e}partition des nombres superabondants'', 
Bull. Soc. Math. France, 103, 1975, 65--90.

\bibitem{ch18:bib4} P. Erd\H{o}s, P. Turan, 
``On some problems of a statistical group theory'', I to VII , 
Zeitschr. fur Wahrschenlichkeitstheorie und verw. Gebiete, 4, 1965, 175--186; 
Acta Math. Hung., 18, 1967, 151-163; 
Acta Math. Hung., 18, 1967, 309--320; 
Acta Math. Hung., 19, 1968, 413--435;
Periodica Math. Hung., 1, 1971, 5-13; 
J. Indian Math. Soc., 34, 1970, 175--192; 
Periodica Math. Hung., 2, 1972, 149--163.

\bibitem{ch18:bib26} J. Grantham, ``The largest prime dividing the
  maximal order of an element of $S_n$'', 
Math.  Comp., 64, 1995, 407--410.

\bibitem{ch18:bib5} E. Landau, 
``Uber die Maximalordung der Permutation gegebenen Grades'',
Handbuch der Lehre von der Verteilung der Primzahlen, vol. 1, 2nd edition, 
Chelsea, New-York, 1953, 222--229.

\bibitem{ch18:bib6} J. P. Massias, 
``Majoration explicite de l'ordre maximum d'un \'{e}l\'{e}ment du
groupe sym\'{e}trique'', 
Ann. Fac. Sci. Toulouse Math., 6, 1984, 269--281.

\bibitem{ch18:bib7} J. P. Massias, J.-L. Nicolas, G. Robin, 
``Evaluation asymptotique de l'ordre
maximum d'un \'{e}l\'{e}ment du groupe sym\'{e}trique'', 
Acta Arithmetica, 50, 1988, 221--242.

\bibitem{ch18:bib8} J. P. Massias, J.-L. Nicolas, G. Robin, 
``Effective bounds for the Maximal Order
of an Element in the Symmetric Group'', 
Math.  Comp., 53, 1989, 665--678.

\bibitem{ch18:bib9} W. Miller, 
``The Maximum Order of an Element of a Finite Symmetric Group'',
Amer. Math. Monthly, 94, 1987, 497--506.

\bibitem{ch18:bib10} H. L. Montgomery, R. C. Vaughan, 
``The large sieve'', 
Mathematika, 20, 1973, 119--134.

\bibitem{ch18:bib11} J.-L. Nicolas, 
``Sur l'ordre maximum d'un \'{e}l\'{e}ment dans le groupe $S_n$ des
permutations'', 
Acta Arithmetica, 14, 1968, 315--332.

\bibitem{ch18:bib12} J.-L. Nicolas, 
``Ordre maximum d'un \'{e}l\'{e}ment du groupe de permutations et
highly composite numbers'', 
Bull. Soc. Math. France, 97, 1969, 129--191.

\bibitem{ch18:bib28} J.-L. Nicolas, 
``Ordre maximal d'un \'{e}l\'{e}ment d'un groupe  de permutations'',  
C.R. Acad. Sci. Paris, 270, 1970, 1473--1476.

\bibitem{ch18:bib13} J.-L. Nicolas, 
``R\'{e}partition des nombres largement compos\'{e}s'', 
Acta Arithmetica, 34, 1979, 379--390.

\bibitem{ch18:bib14} S. Ramanujan, 
``Highly composite numbers'', 
Proc. London Math. Soc., Series 2, 14, 1915, 347--400; 
and ``Collected papers'', Cambridge at the University Press, 1927,
78--128. 

\bibitem{ch18:bib22} S. Ramanujan,  ``Highly composite numbers, 
annotated and with a foreword by J.-L. Nicolas and G.~Robin'', 
Ramanujan J., 1, 1997, 119--153.

\bibitem{ch18:bib15} A. Selberg, ``On the normal density of primes in
  small intervals and the difference between consecutive primes'', 
Arch. Math. Naturvid, 47, 1943, 87--105.

\bibitem{ch18:bib16} S. Shah, 
``An Inequality for the Arithmetical  Function $g(x)$'', 
J. Indian Math. Soc., 3, 1939, 316--318.

\bibitem{ch18:bib17} H. Siebert, 
``Montgomery's weighted sieve for dimension two'', 
Monatsch., Math., 82, 1976, 327--336.

\bibitem{ch18:bib23} N. J. A. Sloane,
``The On-Line Encyclopedia of Integer Sequences'', 
\tt{http://oeis.org}. \rm Accessed 12 December 2012.

\bibitem{ch18:bib25} M. Szalay, 
``On the maximal order in $S_n$ and $S_n^*$'', 
Acta Arithmetica, 37, 1980, 321--331.

\bibitem{ch18:bib24} \tt{http://math.univ-lyon1.fr/$\sim$nicolas/landaug.html}.

\rm 

 

\end{thebibliography}
\end{document}